\documentclass[
11pt]{amsart}

\usepackage{amsthm, amsfonts, amssymb}

\theoremstyle{definition}
\newtheorem{ntn}{Notation}[section]
\newtheorem{dfn}[ntn]{Definition}
\theoremstyle{plain}
\newtheorem{lem}[ntn]{Lemma}
\newtheorem{prp}[ntn]{Proposition}
\newtheorem{thm}[ntn]{Theorem}
\newtheorem{cor}[ntn]{Corollary}
\newtheorem{conj}[ntn]{Conjecture}
\theoremstyle{remark}
\newtheorem{rem}[ntn]{Remark}
\newtheorem{exa}[ntn]{Example}

\def\floor[#1]{\lfloor #1 \rfloor }

\newcommand{\z}{\mathbb{Z}}
\newcommand{\F}{\mathbb{F}}
\newcommand{\pp}{\mathbb{P}}
\newcommand{\q}{\mathbb{Q}}
\newcommand{\R}{\mathbb{R}}
\newcommand{\C}{\mathbb{C}}

\newcommand{\lan}{\langle}
\newcommand{\ran}{\rangle}

\newcommand{\GL}{\mathit{{\rm GL}}}
\newcommand{\SL}{\mathit{{\rm SL}}}

\newcommand{\ppp}{\mathfrak{p}}

\renewcommand{\t}{\mathfrak{t}}

\renewcommand{\u}{\mathcal{U}}

\renewcommand{\top}{{\rm top}}
\newcommand{\inc}{{\rm inc}}
\newcommand{\id}{{\rm id}}

\newcommand{\tors}{{{\rm Tor}_1^{\z}}}

\newcommand{\si}{\sigma}

\newcommand{\del}{\delta}

\newcommand{\arr}{\rightarrow}
\newcommand{\larr}{\longrightarrow}
\newcommand{\harr}{\hookrightarrow}

\newcommand{\se}{\subseteq}

\newcommand{\mt}{\mapsto}
\newcommand{\two}{\twoheadrightarrow}

\newcommand{\Spec}{{\rm Spec}}

\newcommand{\fn}{F^n}
\newcommand{\ffn}{F^{n-2}}
\newcommand{\fff}{{F^\ast}}

\newcommand{\rr}{{F^\ast}}

\newcommand{\stabe}{{\rm Stab}}
\renewcommand{\char}{{\rm char}}
\newcommand{\diag}{{\rm diag}}
\renewcommand{\ker}{{\rm ker}}
\newcommand{\coker}{{\rm coker}}
\newcommand{\im}{{\rm im}}

\newcommand{\Fbar}{\overline{\F}_q}

\newtheoremstyle{athm}
  {}
  {}
  {\itshape}
  {}
  {\scshape}
  {}
  {.5em}
  {\thmnote{#3}}
\theoremstyle{athm}
\newtheorem*{athm}{}

\begin{document}

\title{Homology of $\GL_n$
over algebraically closed fields}
\author{B. Mirzaii}

\begin{abstract}
In this paper we define higher pre-Bloch groups $\ppp_n(F)$ of a
field $F$. When our base field is algebraically closed we study
its connection to the homology of the general linear groups with
finite coefficient $\z/l\z$ where $l$ is a positive integer. As a
result of our investigation we give a necessary and sufficient
condition for the map $H_n(\GL_{n-1}(F), \z/l\z) \arr
H_n(\GL_{n}(F), \z/l\z)$ to be bijective. We prove that this map
is bijective for $n \le 4$. We also demonstrate that the
divisibility of $\ppp_n(\C)$ is equivalent to the validity of the
Friedlander-Milnor Isomorphism Conjecture for $(n+1)$-th homology
of $\GL_n(\C)$.
\end{abstract}

\maketitle

\section*{Introduction}

A theorem of Bloch and Wigner, unpublished and in a somewhat
different form, asserts the existence of the following exact
sequence
\[
0 \arr \q/\z \arr H_3(\SL_2(\C), \z) \arr \ppp(\C) \arr
\bigwedge{}_\z^2
\C^\ast \arr K_2(\C) \arr 0.
\]
A similar exact sequence can be obtained for any algebraically
closed field.
We refer the reader to \cite[Appendix A]{dupont-sah1982} for a proof of
the above exact sequence and for a precise description of the groups and the
maps involved (see also \cite[2.12, 2.14]{sah1989}).

The group $\ppp(\C)$ is called {\it pre-Bloch group} and it has
been the source of
many interesting ideas and connections.
The pre-Bloch group $\ppp(F)$ plays a
very important role in the study of scissors congruences of
polyhedra in connection with Hilbert's third problem
\cite{dupont2001},\cite{sah1989}, it is related very closely to the
third $K$-group $K_3(F)$
which was the main driving force behind Suslin's solution of the
Quillen-Lichtenbaum Conjecture \cite{suslin1991}, and it is used
in establishing certain cases of the Friedlander-Milnor
Isomorphism Conjecture \cite{milnor1983} for certain lower
homology groups \cite{dupont2001}, and so on.

Thus it is natural to ask whether there is a general notion of
higher pre-Bloch groups, and if so, if it carries useful
information. In \cite[Section 4.4]{loday1987} Loday defines a
higher version of the pre-Bloch group, which we denote by
$\ppp_n(F)$, such that $\ppp(F)=\ppp_2(F)$, and predicts that it
should have a close relation with the homology of general linear
groups. We call $\ppp_n(F)$ the $n$-{\it th pre-Bloch group} of
the field $F$.

Although the definition of $\ppp_n(F)$ is easy, which is in terms
of generators and relations, it is difficult to study it directly.
In this article we explore its connection with the homology of the
general linear groups. As we will see this connection is very
close. Here is our main result.

\begin{athm}[{\bf Theorem \ref{inj}.}]
Let $F$ be an algebraically closed field
and let $l$ be a positive integer.
The following conditions are equivalent
\par {\rm (i)} $H_{m}(\GL_{n-1},\z/l\z) {\arr} H_{m}(\GL_{n},\z/l\z)$
is injective for $m=n$ and is surjective for $m=n+1$,
\par {\rm (ii)}
$\ppp_n(F)\otimes \z/l\z=
\begin{cases}
\z/l\z & \text{if $n$ is odd}\\
0 &\text{if $n$ is even.}
\end{cases}$
\end{athm}

It is known by a theorem of Dupont and Sah that for an
algebraically closed field $F$, $\ppp_2(F)$ is divisible
\cite{dupont-sah1982}. The above theorem suggests that
a general version of this fact might be true.

\begin{athm}[{Conjecture \ref{pre-Bloch}.}]
Let $F$ be an algebraically closed field and let $l$ be a
positive integer. Then
\[
\ppp_n(F)\otimes \z/l\z=
\begin{cases}
\z/l\z & \text{if $n$ is odd}\\
0 &\text{if $n$ is even.}
\end{cases}
\]
\end{athm}

From the result of Dupont and Sah the conjecture is
true for n=2 \cite{dupont-sah1982}. In this article
we prove it for $n=3, 4$
and also for all $n$ over the algebraic closure of finite
fields. For the latter we use a result of Friedlander
concerning the homology of general linear groups over algebraic
closure of a finite field.

Here is a strong support for our conjecture.

\begin{athm}[{\bf Proposition \ref{FMC-pn}.}]
Let $l$ be a positive integer. The following conditions are
equivalent
\par {\rm (i)}
$\ppp_n(\C)\otimes \z/l\z=
\begin{cases}
\z/l\z & \text{if $n$ is odd}\\
0 &\text{if $n$ is even}
\end{cases}$
for all $n$,
\par {\rm (ii)} $H_n(\GL_{n-1}(\C), \z/l\z) \overset{\sim}{\larr}
H_n(\GL_{n}(\C), \z/l\z)$ for all $n$,
\par {\rm (iii)}  $H_n(B\GL_{n-1}(\C), \z/l\z) \overset{\sim}{\larr}
H_n( B\GL_{n-1}(\C)^\top, \z/l\z)$ for all $n$.
\end{athm}
~\\
Here $B\GL_{n-1}(\C)^\top$ is the classifying space of
$\GL_{n-1}(\C)$ with it usual topology and $B\GL_{n-1}(\C)$ is the
classifying space of $\GL_{n-1}(\C)$ with $\GL_{n-1}(\C)$ as a
discrete group. The condition (iii) is a special case of the
Friedlander-Milnor Conjecture on the homology of Lie groups with
finite coefficients (see \ref{FMC}).

We briefly outline the organization of the present paper.

In Section \ref{spectral} we introduce a spectral sequence which
will be our main tool in handling the homology of general linear
groups. In this section we will prove an important lemma, which is
used in the proof of Theorem \ref{inj}.

In Section \ref{Pre-Bloch} we define the higher pre-Bloch groups
$\ppp_n(F)$ and give some of its properties. In defining these
groups we follow Suslin's approach for the definition of $\ppp(F)$
in \cite{suslin1991}. Here we also give some satisfactory
description of $\ppp_3(F)$ and $\ppp_4(F)$.

In Section \ref{finite-coeff} we prove Theorem \ref{inj} and
Proposition \ref{FMC-pn}. Here we also prove that Conjecture
\ref{pre-Bloch} is true for algebraic closure of a finite field.


In Section \ref{lower-homology} we show that condition
 (i) or (ii) of Theorem
\ref{inj} is satisfied for $n \le 4$. Here we also establish a new
case of the Friedlander-Milnor conjecture for the fourth homology
of $\GL_3(\C)$ and $\SL_3(\C)$.

In Section \ref{generalization} some of these ideas are
generalized.

\subsection*{Notation}
Here we establish some notations that is used throughout the paper.
In this paper by $H_i(G)$ of a group $G$ we mean the integral homology
group $H_i(G, \z)$. By $\GL_n$ we mean the general linear group $\GL_n(F)$,
where $F$ is an infinite field. By $k$ we mean $\z$ or a prime field and
by $\z/l$ we mean $\z/l\z$, where $l$ is a positive integer.
If $A \arr A'$ is a homomorphism of abelian groups, by $A'/A$ we mean
$\coker(A \arr A')$.

\section{The spectral sequences}\label{spectral}

Let $C_h(\fn)$ be the free $k$-module with a basis consisting of
$(\lan v_0\ran, \dots, \lan v_h\ran)$, where the vectors $v_0,
\dots, v_h \in F^n$ are in general positions, that is every $\min\{h+1,
n\}$ of them is linear independent. By $\lan v_i\ran$ we mean the
line passing through  vectors $v_i$ and $0$.
Let $\partial_0: C_0(\fn) \arr
C_{-1}(\fn):=k$, $\sum_i n_i(\lan v_i\ran) \mt \sum_i n_i$ and
$\partial_h=\sum_{i=0}^h(-1)^id_i: C_h(F^n) \arr C_{h-1}(F^n)$,
$h\ge 1$, where $d_i((\lan v_0\ran, \dots, \lan v_h\ran))= (\lan
v_0 \ran, \dots,\widehat{\lan v_i \ran}, \dots, \lan v_h\ran)$. It
is easy to see that the complex
\begin{equation}\label{exact1}
\hspace{-0.5 cm}
0 \leftarrow k \leftarrow C_0(\fn)
\overset{\partial_1}{\leftarrow} \cdots
\overset{\partial_{n-1}}{\leftarrow} C_{n-1}(\fn)
\overset{\partial_n}
{\leftarrow} C_{n}(\fn)\leftarrow  \cdots
\end{equation}
is exact. Consider the following exact sequence
\begin{equation}\label{exact2}
0 \leftarrow k \leftarrow C_0(\fn)
 \overset{\partial_1}{\leftarrow} \cdots
\overset{\partial_{n-1}}{\leftarrow}
C_{n-1}(\fn) \leftarrow  H_{n-1}(X_n, k) \leftarrow 0,
\end{equation}
where $H_{n-1}(X_n, k):=\ker(\partial_{n-1})$.
We consider $C_{i}(F^n)$ as left
$\GL_n$-module in a natural way.
If it is necessary we convert this action to the right action
by the definition $m.g:=g^{-1}m$ for $g \in \GL_n$ and $m \in C_{i}(F^n)$.

\begin{rem}\label{Xn}
Let $\u(\pp^{n-1})$ be the simplicial set whose for
$0 \le h \le n-1$ non-degenerate
$h$-simplices
are of the form $(\lan v_0\ran, \dots, \lan v_h\ran)$ as in the
above and whose face operators are given by $d_i$. Let $X_n$ be
the geometric realization of $\u(\pp^{n-1})$. It is well-known
that the complex
\[
0 \leftarrow
C_0(\fn)  \overset{\partial_1}{\leftarrow} \cdots
\overset{\partial_{n-1}}{\leftarrow}C_{n-1}(\fn) \leftarrow  0
\]
computes the homology of $X_n$ with coefficient in $k$.  Hence
$H_0(X_n, k)=k$, $H_i(X_n, k)=0$ if $i \neq 0, n-1$ and
$H_{n-1}(X_n, k)=\ker(\partial_{n-1})$.
\end{rem}

The exact sequence
(\ref{exact2})
induces a first quadrant spectral sequence converging to zero with
\begin{gather*}
\begin{array}{l}
%
E_{p, q}^1(n)= \begin{cases} H_q(\rr^p \times \GL_{n-p}, k) &
\text{if $0 \le p \le n$}\\
H_q(\GL_n, H_{n-1}(X_n, k)) & \text{if $p=n+1$} \\
0 & \text{if $p \ge n+2.$}\end{cases}
\end{array}
\end{gather*}
\\ For $1 \le p \le n$, and $q \ge 0$
the differential
$d_{p, q}^1(n)$
equals $\sum_{i=1}^p(-1)^{i+1}H_q(\alpha_{i, p})$ where
$\alpha_{i, p}: \fff^p \times \GL_{n-p} \arr \fff^{p-1} \times \GL_{n-p+1}$,
\[
(a_1, \dots, a_p, A) \mt
(a_1, \dots, \widehat{a_i}, \dots, a_p,
\left(
\begin{array}{cc}
a_i & 0          \\
0   & A
\end{array}
\right))
\]
(see the proof of \cite[Thm. 3.5]{mirzaii2003} for details).
In particular for $0 \le p \le n$,
\[
d_{p, 0}^1(n)= \begin{cases} {\rm id}_k &
\text{if $p$ is odd}\\
0 & \text{if $p$ is even,}\end{cases}
\]
so $E_{p, 0}^2(n)=0$ for $ p \neq n, n+1$.
In fact this is also true for $p=n, n+1$.
Applying the right exact functor $H_0$ to the exact sequence
\[
C_{n+1}(\fn) \arr C_{n}(\fn) \arr H_{n-1}(X_n, k) \arr 0
\]
we get the exact sequence
\[
H_0(\GL_n, C_{n+1}(\fn)) \arr H_0(\GL_n, C_{n}(\fn)) \arr
H_0(\GL_n, H_{n-1}(X_n, k)) \arr 0.
\]
The group $\GL_n$ acts transitively on the basis
$(\lan v_0\ran, \dots, \lan v_n\ran)$ of $C_{n}(\fn)$ so
$H_0(\GL_n, C_{n}(\fn))=k$. From this we obtain
\begin{equation}\label{h1gn}
H_0(\GL_n, H_{n-1}(X_n,k))=
\begin{cases} 0 &
\text{if $n$ is odd}\\
k & \text{if $n$ is even.}
\end{cases}
\end{equation}

Consider $\ffn$ as a vector subspace of $\fn$ generated
by $e_3, e_4, \dots, e_{n}$ (so $\GL_{n-2}$ embeds in $\GL_n$ as
$\diag(1, 1, \GL_{n-2})$). Let
${L'}_\ast$ and $L_\ast$ be the complexes
\begin{gather*}
0  \leftarrow  0  \leftarrow   \ \ \  0  \ \ \    \leftarrow  \ \ \ \  k
\ \ \
 \leftarrow  C_0(\ffn)  \leftarrow
\cdots \leftarrow  H_{n-3}(X_{n-2},k)  \leftarrow   0
\\
0 \leftarrow k \leftarrow C_0(\fn) \leftarrow C_1(\fn)
 \leftarrow  C_2(\fn) \leftarrow
\cdots \leftarrow  H_{n-1}(X_{n},k) \leftarrow   0
\end{gather*}
respectively, that is $L_0'=0$, $L_1'=0$, $L_2'=k$, $L_{i+3}'=C_i(\ffn)$
for $i=0, \dots, n-3$, $L_{n+1}'=H_{n-3}(X_{n-2}, k)$ and $L_{i}'=0$
for $i \ge n+2$, $L_0=k$, $L_{i+1}=C_i(\fn)$ for $i=0, \dots, n-1$,
$L_{n+1}=H_{n-1}(X_{n}, k)$ and $L_{i}=0$
for $i \ge n+2$. Define the map of complexes
${L}_\ast' \overset{\theta_\ast}{\arr} {L}_\ast$, given by
\begin{gather*}
(\lan v_1\ran, \dots, \lan v_j\ran) \overset{\theta_j}{\mt}
(\lan e_1\ran,\lan e_2\ran,\lan v_1\ran, \dots, \lan v_j\ran)-
(\lan e_1\ran, \lan e_1+e_2\ran,\lan v_1\ran, \dots, \lan v_j\ran)\\
\!\!\!\!\!\!\!\!\!\!\!\!\!\!
+(\lan e_2\ran, \lan e_1+e_2\ran, \lan v_1\ran, \dots, \lan v_j\ran).
\end{gather*}
This induces a map of bicomplexes
\begin{gather*}
{L}_\ast' \otimes_{\GL_{n-2}} {F}_\ast' \arr
{L}_\ast \otimes_{\GL_{n}} {F}_\ast \arr
{L}_\ast \otimes_{\GL_{n}} {F}_\ast /{L}_\ast'
\otimes_{\GL_{n-2}} {F}_\ast',
\end{gather*}
where $F_\ast: \cdots \arr F_2 \overset{\delta_2}{\arr}
F_1 \overset{\delta_1}{\arr} F_0 \arr 0$ is a
free left $\GL_n$-resolution of $k$
and ${F}_\ast'$ is ${F}_\ast$ as $\GL_{n-2}$-resolution.
Thus one gets the map of spectral sequences
\begin{gather*}
{E'}_{p, q}^r(n) \arr {E}_{p, q}^r(n) \arr {E''}_{p, q}^r(n),
\end{gather*}
where all the three spectral sequences converge to zero.
By a similar approach as in the proof of \cite[Thm. 3.5]{mirzaii2003},
one sees that the spectral sequence ${E'}_{p, q}^1(n)$ is of the form
\begin{gather*}
{E'}_{p, q}^1(n)= \begin{cases}
E_{p-2, q}^1(n-2) &
\text{if $p \ge 2$}\\
0 & \text{if $p=0, 1.$} \end{cases}
\end{gather*}

It is not difficult to see that for $2 \le p \le n$ and $q \ge 0$ the map
${E'}_{p, q}^1(n) \arr {E}_{p, q}^1(n)$ is the map $H_q(\inc)$  induced by
$\inc:\fff^{p-2} \times \GL_{n-p} \arr \fff^{p} \times \GL_{n-p}$ with
$A \mt (1, 1 , A)$ and by a little work one sees that
\begin{gather*}
{E''}_{p, q}^1(n)={E}_{p, q}^1(n)/{E'}_{p, q}^1(n)
\end{gather*}
for $0 \le p \le n$ (see \cite[Section 4]{mirzaii2003}).

From the exact sequence of complexes
\begin{gather*}
0 \arr {L}_p' \otimes_{\GL_{n-2}} {F}_{\ast}' \arr
{L}_p \otimes_{\GL_{n}} {F}_\ast \arr
{L}_p \otimes_{\GL_{n}} {F}_\ast /{L}_p \otimes_{\GL_{n-2}}
{F}_\ast'\arr 0
\end{gather*}
we obtain the long exact sequence
\begin{equation}\label{exact4}
\hspace{-2.5 cm}
\cdots \arr  {E'}_{p, q}^1(n)
\arr {E}_{p, q}^1(n)
\arr {E''}_{p, q}^1(n)
\arr {E'}_{p, q-1}^1(n)
\end{equation}
\vspace{-0.6 cm}
\begin{gather*}
\vspace{-0.3 cm}
\hspace{4 cm}
\arr \cdots
\arr {E'}_{p, 0}^1(n)
\arr {E}_{p, 0}^1(n)
\arr {E''}_{p, 0}^1(n)
\arr 0.
\end{gather*}
This exact sequence is studied in the above if $0 \le p \le n, q
\ge 0$ and $p=n+1, q=0$. We will come back to it later.

Here is an important lemma which is used in the proof of Theorem \ref{inj}.

\begin{lem}\label{h1g}
$E_{n, i}^2(n)=0$ for all $i$. In particular
${E'}_{n, i}^2(n)={E''}_{n, i}^2(n)=0$.
\end{lem}
\begin{proof}
 Consider the following commutative diagram with exact columns
\begin{gather*}
\begin{array}{ccccc}
0              & &     0          & &      0          \\
\Big\downarrow & & \Big\downarrow & & \Big\downarrow  \\
H_{n-1}(X_n, k)\otimes_{\GL_n}F_{i+1} & \!\!\!\!\!\arr &
\!\!\!\!\! H_{n-1}(X_n, k) \otimes_{\GL_n}F_i&
\!\!\!\!\!\arr &  \!\!\!\!\! H_{n-1}(X_n, k)\otimes_{\GL_n}F_{i-1}  \\
\Big\downarrow & & \Big\downarrow & & \Big\downarrow  \\
C_{n-1}(\fn)\otimes_{\GL_n} F_{i+1} & \!\!\!\!\!\arr &
\!\!\!\!\!C_{n-1}(\fn)\otimes_{\GL_n}F_i &
\!\!\!\!\!\arr &  \!\!\!\!\! C_{n-1}(\fn)\otimes_{\GL_n}F_{i-1}  \\
\Big\downarrow & & \Big\downarrow & & \Big\downarrow  \\
C_{n-2}(\fn)\otimes_{\GL_n}F_{i+1}  & \!\!\!\!\!\arr &
\!\!\!\!\! C_{n-2}(\fn) \otimes_{\GL_n}F_i&
\!\!\!\!\!\arr & \!\!\!\!\! C_{n-2}(\fn) \otimes_{\GL_n}F_{i-1}.\\
&&&&
\end{array}
\end{gather*}
Set $\si_{j}=(\lan e_1 \ran, \dots, \lan e_{j+1} \ran)$ for $0 \le
j \le n-1$. Let $\si_{n-1}\otimes x  \in C_{n-1}(\fn)
\otimes_{\GL_n}F_i$ represent an element of the group
$H_i(\GL_n,C_{n-1}(\fn)) \simeq  H_i(\fff^n, k)$ such that $d_{n,
i}^1(\overline{\si_{n-1}\otimes x})=0$. Then
\begin{gather*}
(\partial_{n-1}\otimes \id_{F_i})(\si_{n-1} \otimes x)=
\partial_{n-1}(\si_{n-1})\otimes x
\in \im( \id_{C_{n-2}}\otimes \del_{i+1}).
\end{gather*}
Let $ \si_{n-2}\otimes \del_{i+1}(y)=\partial_{n-1}(\si_{n-1})\otimes x$.
It is easy to see that
\[
\si_{n-2} \otimes \del_{i+1}(y)=
\si_{n-2}\otimes (\sum_{i=0}^{n-1}(-1)^{i}g_ix),
\]
where $g_i \in \GL_n$ is the permutation matrix such that
\[
g_{i}^{-1}(e_1, \dots, \widehat{e_{i+1}}, \dots, e_n, e_{i+1})=
(e_1, \dots, e_{n}).
\]
The inclusions $\fff^n \se \stabe_{\GL_n}(\si_{n-2}) \se \GL_n$
induce the commutative diagram
\begin{gather*}
\begin{array}{ccccc}
k\otimes_{\fff^n}F_{i+1} & \larr &
 k \otimes_{\stabe_{\GL_n}(\si_{n-2})}F_{i+1}  &
\overset{\sim}{\larr}
& C_{n-2}(\fn) \otimes_{\GL_n} F_{i+1}  \\
\Big\downarrow & & \Big\downarrow & & \Big\downarrow  \\
k\otimes_{\fff^n} F_{i} &
{\larr} &
 k\otimes_{\stabe_{\GL_n}(\si_{n-2}) } F_{i} &
\overset{\sim}{\larr} & C_{n-2}(\fn) \otimes_{\GL_n} F_{i}  \\
\Big\downarrow & & \Big\downarrow & & \Big\downarrow  \\
k\otimes_{\fff^n} F_{i-1} & \larr &
k  \otimes_{\stabe_{\GL_n}(\si_{n-2}) } F_{i-1}  &
\overset{\sim}{\larr} & C_{n-2}(\fn) \otimes_{\GL_n}F_{i-1}.
\end{array}
\end{gather*}
Since $H_i(\fff^n, k) \simeq H_i(\stabe_{\GL_n}(\si_{n-2}), k)$,
there is a $y' \in F_{i+1}$ such that
\[
1 \otimes \del_{i+1}(y') = 1\otimes(\sum_{i=0}^{n-1}(-1)^{i}g_ix)
\in k \otimes_{\fff^n} F_{i}.
\]
Let $\si_{n}=(\lan e_1 \ran, \dots, \lan e_n \ran, \lan e_1+
\cdots+ e_n \ran )$. Clearly
$\partial_{n}(\si_{n}) \in H_{n-1}(X_n, k)$. If
$z=\partial_{n}(\si_{n})\otimes x$, then
\begin{gather*}
(j \otimes \id_{F_i})(z)
=(-1)^{n} \si_{n-1} \otimes x +
\tilde{\si}_{n-1}\otimes (\sum_{i=0}^{n-1}(-1)^{i}g_i x) \\
\hspace{2.4 cm}
=(-1)^{n} \si_{n-1} \otimes x +
\si_{n-1}\otimes (\sum_{i=0}^{n-1}(-1)^{i} gg_i x),
\end{gather*}
where $j:H_{n-1}(X_n, k) \harr C_{n-1}(\fn)$,
$\tilde{\si}_{n-1}=(\lan e_1 \ran, \dots, \lan e_{n-1} \ran,
\lan e_1+\cdots+e_n \ran)$ and $g \in \stabe_{\GL_n}(\si_{n-2})$ with
$g^{-1}\tilde{\si}_{n-1}=\si_{n-1}$.
Since $1\otimes g_i x= 1 \otimes gg_i x$ in
$ k\otimes_{\stabe_{\GL_n}(\si_{n-2}) } F_{i}$, there exist
$y'' \in F_{i+1}$ such that
\[
1\otimes \sum_{i=0}^{n-1}(-1)^{i}gg_i x=
1\otimes \sum_{i=0}^{n-1}(-1)^{i} g_i x + 1 \otimes \delta_{i+1}(y'')
\in k \otimes_{\fff^n} F_i.
\]
Now in $ C_{n-1}(\fn)\otimes_{\GL_n} F_i$ we have
\begin{gather*}
\hspace{-3 cm}
(j\otimes \id_{F_i})(z)=(-1)^{n}\si_{n-1}\otimes x +
\si_{n-1} \otimes \sum_{i=0}^{n-1}(-1)^{i} gg_i x \\
\hspace{2.3 cm}
=(-1)^{n}\si_{n-1} \otimes x+
\si_{n-1} \otimes \sum_{i=0}^{n-1}(-1)^{i} g_i x +
\si_{n-1} \otimes \delta_{i+1}(y'')\\
\hspace{-0.8  cm}
=(-1)^{n} \si_{n-1}\otimes x+ \si_{n-1}\otimes \delta_{i+1}(y'+y'').
\end{gather*}
This completes the proof of the triviality of ${E}_{n, i}^2(n)$.
The triviality of ${E'}_{n, i}^2(n)$ follows immediately from this, because
${E'}_{n, i}^2(n)={E}_{n-2, i}^2(n-2)$. The triviality of
${E''}_{n, i}^2(n)$ follows from these and applying the Snake lemma
to the following commutative diagram with exact rows
\begin{gather*}
\begin{array}{ccccccccc}
& & {E'}_{n+1, i}^1(n) & \arr & {E}_{n+1, i}^1(n)&
\arr & {E''}_{n+1, i}^1(n) & \arr & 0  \\
& &\Big\downarrow & & \Big\downarrow & & \Big\downarrow & & \\
0 & \arr & \ker({d'}_{n,i}^1(n)) &
\arr & \ker ({d}_{n,i}^1(n))  & \arr & \ker ({d''}_{n,i}^1(n))& \arr & 0.
\end{array}
\end{gather*}
\end{proof}

\section{Higher pre-Bloch groups}\label{Pre-Bloch}

In this section we define the higher pre-Bloch groups
$\ppp_n(F)$ and investigate some of its properties.
\begin{dfn}
Set $\t_n^{(k)}(F):= H_1(\GL_{n}, H_{n-1}(X_{n}, k))$. We denote
$\t_n^{(\z)}(F)$ by $\t_n(F)$. By convention $\t_n^{(k)}(F)=0$ for $n=0,1$.
\end{dfn}
From the short exact sequence
\[
0\arr \partial_{n+1}(C_{n+1}(\fn)) \overset{\alpha}{\arr} C_n(\fn)
\overset{\beta}
{\arr} H_{n-1}(X_n, \z) \arr 0
\]
one obtains the long exact sequence
\begin{gather*}
\hspace{-2.5 cm}
\cdots \arr H_1(\GL_n, \partial_{n+1}(C_{n+1}(\fn)))
\overset{H_1(\alpha)}{\larr} H_1(\GL_n,C_n(\fn)) \\
\hspace{4 cm}
\overset{H_1(\beta)}{\larr}  \t_n(F) \arr
H_0(\GL_n, \partial_{n+1}(C_{n+1}(\fn))) \arr \cdots.
\end{gather*}
If $n$ is even, the composition
\[
\fff=H_1(\GL_n,C_n(\fn)) \overset{H_1(\beta)}{\arr} \t_n(F)
\overset{H_1(j)}{\larr}
H_1(\GL_n,C_{n-1}(\fn))=\fff^n
\]
is given by $H_1(j)\circ H_1(\beta)(a)=(a,\dots, a) \in \fff^n$,
where $j:H_{n-1}(X_n, k) \harr C_{n-1}(\fn)$. Thus
$H_1(\beta)$ is injective. If $n$ is odd, the composition
\[
H_1(\GL_n,C_{n+1}(\fn)) 
{\arr} H_1(\GL_n, \partial_{n+1}(C_{n+1}(\fn)))
\overset{H_1(\alpha)}{\larr}
H_1(\GL_n,C_{n}(\fn))
\]
is surjective and so $H_1(\alpha)$ is surjective.
Now from the
above long exact sequence one gets the following exact sequences
\begin{gather*}
0 \arr \fff \arr \t_n(F) \arr H_0(\GL_n, \partial_{n+1}(C_{n+1}(\fn)))
\arr 0, \ \
{\rm if\ {\it n}\ is\  even},
\end{gather*}
\begin{gather*}
0 \arr \t_n(F) \arr H_0(\GL_n, \partial_{n+1}(C_{n+1}(\fn))) \arr \z \arr 0,
\ \ \
{\rm if\ {\it n}\ is \ odd}.
\end{gather*}
To study $H_0(\GL_n, \partial_{n+1}(C_{n+1}(\fn)))$
apply the functor $H_0$ to
\[
C_{n+2}(\fn)\arr C_{n+1}(\fn) \arr \partial_{n+1}(C_{n+1}(\fn)) \arr 0.
\]
Thus we get the exact sequence
\[
(C_{n+2}(\fn))_{\GL_n} \arr (C_{n+1}(\fn))_{\GL_n} \arr
H_0(\GL_n, \partial_{n+1}(C_{n+1}(\fn))) \arr 0.
\]
Let $E_n=\sum_{i=1}^ne_i \in \fff^n$, $a=\sum_{i=1}^na_ie_i\in \fff^n$, where
$a_i \in \fff- \{1 \}$ and $a_i \neq a_j$ if $i \neq j$.
Denote the orbit of the frame
$(\lan e_1 \ran, \dots, \lan e_n  \ran, \lan E_n  \ran,\lan a  \ran )
\in C_{n+1}(\fn)$
by $p(a)$ and orbit of the frame
$(\lan e_1 \ran, \dots, \lan e_n  \ran, \lan E_n  \ran,\lan a  \ran,
\lan b  \ran )\in C_{n+2}(\fn)$ by $p(a, b)$,
where $b=\sum_{i=1}^nb_ie_i\in \fff^n$,
$b_i \in \fff- \{1 \}$, $b_i \neq b_j$ if $i \neq j$ and
$a_i \neq b_j$ for all $i, j$. We see that
\begin{gather*}
(C_{n+1}(\fn))_{\GL_n}=\coprod_{a} \z.p(a), \ \ \
(C_{n+2}(\fn))_{\GL_n}=\coprod_{a, b} \z.p(a, b).
\end{gather*}
A direct computation shows that
\begin{gather*}
\overline{\partial_{n+2}}(p(a, b))= \sum_{i=1}^n (-1)^{i+1}
p((\frac{b_1-b_i}{a_1-a_i}, \dots, \widehat{\frac{b_i-b_i}{a_i-a_i}},
\dots, \frac{b_n-b_i}{a_n-a_i},\frac{b_i}{a_i}))\\
\hspace{1.5 cm}
+ (-1)^n p((\frac{b_1}{a_1}, \dots, \frac{b_n}{a_n} ))-
(-1)^n p(a)+(-1)^n p(b).
\end{gather*}
Therefore $H_0(\GL_n, \partial_{n+1}(C_{n+1}(\fn)))$ is generated by
$[a_1; \dots; a_n] \in \pp^{n-1}$, $a_i \in \fff- \{1 \}$,
$a_i \neq a_j$ if $i \neq j$, and relations
\begin{gather*}
\hspace{-5 cm}
[b_1; \dots; b_n]-[a_1; \dots; a_n]+
[\frac{b_1}{a_1}; \dots; \frac{b_n}{a_n}] \\
\hspace{2 cm}
-\sum_{i=1}^n (-1)^{i+n}[
\frac{b_1-b_i}{a_1-a_i}; \dots; \widehat{\frac{b_i-b_i}{a_i-a_i}};
\dots; \frac{b_n-b_i}{a_n-a_i};\frac{b_i}{a_i}]=0,
\end{gather*}
where $b_i$ are as above. If in the above we replace $a_i/a_n$ and
$b_i/b_n$ with $a_i$ and $b_i$ respectively, one sees that the
group $H_0(\GL_n, \partial_{n+1}(C_{n+1}(\fn)))$ is generated by
the symbols $[a_1, \dots, a_{n-1}]$, $a_i \in \fff- \{1 \}$, $a_i
\neq a_j$ if $i \neq j$ and relations
\begin{gather*}
[b_1, \dots, b_{n-1}]-[a_1, \dots, a_{n-1}]+
[\frac{b_1}{a_1}, \dots, \frac{b_{n-1}}{a_{n-1}}]-
[\frac{b_1-1}{a_1-1}, \dots,
\frac{b_{n-1}-1}{a_{n-1}-1}]\\
-\sum_{i=1}^{n-1} (-1)^{i+n}[
\frac{b_1b_i^{-1}-1}{a_1a_i^{-1}-1}, \dots,
\widehat{\frac{b_ib_i^{-1}-1}{a_ia_i^{-1}-1}},
\dots, \frac{b_{n-1}b_i^{-1}-1}{a_{n-1}a_i^{-1}-1},
\frac{b_i^{-1}-1}{a_i^{-1}-1}]=0,
\end{gather*}
where $b_i \in \fff- \{1 \}$, $b_i \neq b_j$ if $i \neq j$ and
$a_i \neq b_j$ for all $i, j$. When $n=2$ this is the definition
of $\ppp(F)$.

Thus one can think of
$H_0(\GL_n, \partial_{n+1}(C_{n+1}(\fn)))$ as a natural generalization
of $\ppp(F)$ for $n \ge 3$. So we allow ourself to make the following
definition.

\begin{dfn}
The group $H_0(\GL_n, \partial_{n+1}(C_{n+1}(\fn)))$ is called the $n$-th
pre-Bloch group of $F$ and we denote it by $\ppp_n(F)$.
\end{dfn}

{}From the above
we have the following exact sequences
\begin{equation}\label{even}
0 \arr \fff \arr \t_n(F) \arr \ppp_n(F) \arr 0, \ \
{\rm if\ {\it n}\ is\  even},\\
\end{equation}
\begin{equation}\label{odd}
0 \arr \t_n(F) \arr \ppp_n(F) \arr \z \arr 0,
\ \ \
{\rm if\ {\it n}\ is \ odd}.
\end{equation}

\begin{rem}\label{pn-phin}
In \cite[2.7]{yagunov2000} Yagunov defines another version of
higher pre-Bloch groups, denoted by $\wp^n(F)$. He also defines
the classical pre-Bloch group ${\wp^n(F)}_{\rm cl}$. Our
definition of the pre-Bloch group is very close to his definition
of the classical pre-Bloch group. In fact
\[
{\wp^n(F)}_{\rm cl}=
\begin{cases}
\ker(\ppp_n(F) \arr \z) & \text{if $n$ is odd}\\
\ppp_n(F)&\text{if $n$ is even.}
\end{cases}
\]
See \cite[3.11]{yagunov2000} for the relation between
${\wp^n(F)}_{\rm cl}$ and $\wp^n(F)$.
\end{rem}

Since ${E}_{n+1, 1}^1(n)=\t_n(F)$ for $k=\z$, from the exact
sequence (\ref{exact4}) we have the following exact sequence
\begin{equation}\label{exact5}
\t_{n-2}(F) \arr \t_n(F) \arr {E''}_{n+1, 1}^1(n) \arr 0.
\end{equation}
An easy calculation shows that $\ker(d_{n, 1}^1(n)) \se \fff^n$
is generated by elements of the form
\[
A= \begin{cases}
(a_1, a_1, \dots, a_j, a_j)
&\text{if $n=2j$ }\\
(a_1, 1, \dots, a_{j}, 1, \prod_{i=1}^j a_{i}^{-1})
& \text{if $n=2j+1.$ }\end{cases}
\]
This proves that
$\ker(d_{n, 1}^1(n))\simeq \fff^{\floor [n/2]}$. So we
have a surjective map $\t_n(F) \arr \fff^{\floor [n/2]}$.
Using (\ref{exact5}) we obtain a surjective map
${E''}_{n+1, 1}^1(n) \arr \fff$. It is not difficult to see
that the composition $\t_n(F) \arr {E''}_{n+1, 1}^1(n) \arr \fff$
splits the exact sequence (\ref{even}) for $n$ even.
So we have proved the following lemma.

\begin{lem}\label{tn-pn}
Let $n \ge 2$. $\t_n(F) \simeq \fff \oplus \ppp_n(F)$ if $n$ is even and
$\ppp_n(F) \simeq \z \oplus \t_n(F)$ if $n$ is odd.
\end{lem}

In the following lemma we give some satisfactory description of
the group $\t_n(F)$, $n=3,4$, for arbitrary infinite field $F$.
This also gives a better description of
$\ppp_n(F)$ for $n= 3, 4$.

\begin{lem}\label{p234}
Let $F$ be an infinite field. Then
\par {\rm (i)} $\t_2(F)\simeq \fff \oplus \ppp_2(F)$,
\par {\rm (ii)} $\t_3(F) \simeq \fff$, therefore
$\ppp_3(F) \simeq \fff \oplus \z$,
\par {\rm (iii)} $\t_4(F) \simeq \fff \oplus \ppp_4(F)$ and
there is an exact sequence
\[
\t_2(F)  \arr \t_4(F) \arr \fff \arr 1.
\]
\end{lem}
\begin{proof}
(i) This part has already been proven in Lemma \ref{tn-pn}. \\
(ii) By Lemma \ref{h1g} and \cite[Cor. 3.5]{mirzaii2004}
 the ${E}_{p, q}^2(3)$-terms are of the following form
\begin{gather*}
\begin{array}{ccccccc}
\ast & \ast    &      &     &              &      &    \\
0    & 0       & \ast &\ast & \ast         &      &    \\
0    & 0       &  0   &  0  & \ast         & 0    &    \\
0    & 0       &  0   &  0  & E_{4,1}^2(3) & 0    &    \\
0    & 0       &  0   &  0  & 0            & 0    &  \cdots.
\end{array}
\end{gather*}
By a similar arguments as in the proof of
\cite[Lemma 3.6]{mirzaii2004} we have $E_{0,4}^3(3)=0$.
Since the spectral sequence converges to zero,
$E_{4, 1}^2(3)=E_{4, 1}^\infty(3)=0$. Therefore $\t_3(F) \simeq \fff$.\\
(iii) For $n=4$ we look at the spectral sequence ${E''}_{p, q}^2(4)$.
By Lemma \ref{h1g} and \cite[Thm. 5.5]{mirzaii2005}
the ${E''}_{p, q}^2(4)$-terms are of the following form
\begin{gather*}
\begin{array}{cccccccc}
\ast & \ast    &      &       &      &                   &      &    \\
0    & 0       & \ast & \ast  & 0    &                   &      &     \\
0    & 0       & \ast & \ast  & 0    & \ast              &      &    \\
0    & 0       &  0   & 0     & 0    & \ast              & 0    &   \\
0    & 0       &  0   & 0     & 0    & {E''}_{5, 1}^2(4) & 0    &    \\
0    & 0       &  0   & 0     & 0    & 0                 & 0    &  \cdots.
\end{array}
\end{gather*}
With a similar argument to the case $n=3$, using the results of
\cite{mirzaii2005} one can show that ${E''}_{0, 5}^3(4)=0$ (see
the proof of \cite[Lemma 3.6]{mirzaii2004}). With a little bit
work one can prove that ${E''}_{2, 3}^2=0$ (see \cite{mirzaii2004}
or \cite{mirzaii2005} to get an idea how one can do that). An easy
analysis of the spectral sequence shows that ${E''}_{5,
1}^2(4)={E''}_{5, 1}^\infty(4)=0$. The exact sequence follows from
this and the exact sequence (\ref{exact5}).
\end{proof}

\section{Homology of $\GL_n$ with finite coefficient}\label{finite-coeff}

In this section we show that over an algebraically
closed field $F$ the pre-Bloch group $\ppp_n(F)$ is closely related to the
homology of $\GL_n$ with finite coefficient.

\begin{lem}\label{m1}
Let $F$ be an infinite field, $k$ a prime field and assume that
$n \ge 3$, $j \ge 0$ be integers such that $n+1 \ge j$.
Let  $H_q(\inc): H_q(\GL_{n-2}, k) \arr H_q(\GL_{n-1}, k)$
be surjective for $0 \le q \le j-1$. Then the following conditions are
equivalent;
\par {\rm (i)}  $H_j(\inc): H_j(\GL_{n-1}, k) \arr H_j(\GL_{n}, k)$
is surjective,
\par {\rm (ii)} $H_j(\inc): H_j(\fff \times \GL_{n-1}, k)
\arr H_j(\GL_{n}, k)$ is surjective.
\end{lem}
\begin{proof}
See \cite[Lem. 4.1]{mirzaii2003}.
\end{proof}

\begin{lem}\label{m2}
Let $F$ be an infinite field, $k$ a prime field and assume that
$n \ge 3$, $j \ge 0$ be integers such that $n \ge j$.
Let $H_q(\inc ): H_q(\GL_{m-1}, k) \arr H_q(\GL_{m}, k)$ be isomorphism
for $m=n, n-1$ and $0 \le q \le {\min}\{j-1, m-2\}$.
Then the following conditions are equivalent;
\par {\rm (i)} $H_j(\inc): H_j(\GL_{n-1}, k) \arr H_j(\GL_{n}, k)$
is bijective,
\par {\rm (ii)} $H_j(\fff^2\times \GL_{n-2}, k) \overset{\tau_2}{\arr}
H_j(\fff \times \GL_{n-1}, k) \overset{\tau_1}{\arr} H_j(\GL_{n}, k) \arr 0$
is exact, where $\tau_1=H_j(\inc)$ and
$\tau_2= H_j(\alpha)- H_j(\inc)$, $\alpha: (a, b, A) \mt (b, \diag(a, A))$.
\end{lem}
\begin{proof}
See \cite[Lem. 4.2]{mirzaii2003}.
\end{proof}

\begin{prp}[Stability]\label{stability}
Let $F$ be an algebraically closed field.
Then $H_{m}(\GL_{n-1},\z/l) {\arr} H_{m}(\GL_{n},\z/l)$
is surjective if $m \le n$  and is injective if $ m \le n-1$.
\end{prp}
\begin{proof}
These results are already known and immediately follow from
Suslin's homological stability theorem \cite[Thm.
3.4]{suslin1985}. But for this special case we give a proof that
is much easier than Suslin's proof. Here we may assume that $l$ is
a prime. The proof is by
induction on $n$. If $n=1$ then everything is obvious. Assume the
induction hypothesis, that is $H_m(\GL_{j-1},\z/l) \arr
H_m(\GL_{j},\z/l)$ is surjective if $m \le j$ and is bijective if
$m \le j-1$, where $1 \le j \le n-1$. Consider the spectral
sequence ${E''}_{p, q}^2(n)$ with $k=\z/l$. It is sufficient to
prove that ${E''}_{p, q}^2(n)=0$ if $p+q \le n+1$, $0 \le q \le
n-2$ and ${E''}_{2, n-1}^2(n)=0$. Because then we obtain
${E''}_{0, m}^2(n)=0$ for $ 0 \le m \le n$ and ${E''}_{1,
m}^2(n)=0$ for $ 0 \le m \le n-1$ and by applying Lemmas \ref{m1}
and \ref{m2} we get the desired results. The proof is analogue
(and even easier) than the proof of \cite[Thm. 4.3]{mirzaii2003}.
So we refer the reader to the proof of that theorem.
Note that here one must use the fact
that $H_{2i+1}(\fff, \z/l)=0$ for $i \ge 0$, since $\fff$ is a
divisible group \cite[Prop. 4.7]{dupont2001}.
\end{proof}

\begin{thm}\label{inj}
Let $F$ be algebraically closed.
The following conditions are equivalent
\par {\rm (i)} $H_{m}(\GL_{n-1},\z/l) {\arr} H_{m}(\GL_{n},\z/l)$
is injective for $m=n$ and is surjective for $m=n+1$,
\par {\rm (ii)}
$\ppp_n(F)\otimes \z/l=
\begin{cases}
\z/l & \text{if $n$ is odd}\\
0 & \text{if $n$ is even.}
\end{cases}$
\end{thm}
\begin{proof}
We may assume that $l$ is a prime. Again we look at the spectral
sequence ${E''}_{p, q}^1(n)$ with $k=\z/l$. By Lemma \ref{h1g},
${E''}_{n, 2}^2(n)=0$. By Proposition \ref{stability} and a
similar argument as the proof of \cite[Thm. 4.3]{mirzaii2003} one
can show that ${E''}_{p, q}^2(n)=0$ for $p+q=n+2$, where $3 \le q
\le n$. Since the spectral sequence converges to zero one sees
that ${E''}_{n+1, 1}^2(n)=0$ if and only if ${E''}_{1, n}^2(n)=0$
and ${E''}_{0, n+1}^2(n)=0$. Note that by Lemma \ref{tn-pn},
$\ppp_n(F)$ has the desired property if and only if
$\t_n(F)\otimes \z/l=0$ if and only if $\t_n^{(\z/l)}(F)=0$ (use
Remark \ref{Xn}). By Proposition \ref{gli} this theorem is true
for $n=2,3$. Thus by induction we may assume that the theorem is
true for lower cases.

By (\ref{exact5}) and the induction step
${E''}_{n+1, 1}^2(n)=0$ if and only if $\t_n^{(\z/l)}(F)=0$.
By Lemmas \ref{m1} and \ref{m2} we have
${E''}_{1, n}^2(n)={E''}_{0, n+1}^2(n)=0$ if and only if
the map $H_{m}(\GL_{n-1},\z/l) {\arr} H_{m}(\GL_{n},\z/l)$
is injective for $m=n$ and is surjective for $m=n+1$.
This completes the proof of the theorem.
\end{proof}

So it is convenient to make the following conjecture (which
easily follows from Conjecture \ref{H(X_N)} using Lemma \ref{tn-pn}).

\begin{conj}\label{pre-Bloch}
Let $F$ be algebraically closed.
Then
\[
\ppp_n(F)\otimes \z/l=
\begin{cases}
\z/l & \text{if $n$ is odd}\\
0 & \text{if $n$ is even.}
\end{cases}
\]
\end{conj}

\begin{rem}
(i) By a theorem of Dupont and Sah Conjecture \ref{pre-Bloch} is
true for $n=2$ \cite[App. A]{dupont-sah1982} and by lemma \ref{p234}
it is also true for $n=3, 4$.
\par (ii) Conjecture \ref{pre-Bloch} is related very closely to a
conjecture of Yagunov \cite[Conj. 0.2]{yagunov2000}. In
\cite{yagunov2000} he defines certain pre-Bloch groups
$\varphi^n(F)$ and conjectures that they are divisible. By Remark
\ref{pn-phin} and \cite[Prop. 3.11]{yagunov2000}, up to
$2$-torsion, Conjecture \ref{pre-Bloch} implies Conjecture 0.2
from \cite{yagunov2000}.
\end{rem}

In the rest of this section we prove certain results that support
Conjecture \ref{pre-Bloch}. First a theorem due to Friedlander.

\begin{thm}\label{friedlander-thm}
Let
$\Fbar$ be the algebraic closure of
the finite field ${\F_q}$. Then
\par {\rm (i)} $H_{i}(\GL_{n}(\Fbar)) \arr  H_{i}(\GL_{n+1}(\Fbar))$ is
isomorphism for $i \le 2n-1$,
\par {\rm (ii)} $H_{i}(\SL_{n}(\Fbar)) \arr  H_{i}(\SL_{n+1}(\Fbar))$ is
isomorphism for $i \le 2n-1$.
\end{thm}
\begin{proof}
See \cite[Thm. 3]{friedlander1976}.
\end{proof}

\begin{cor}\label{pn-fbar}
{\rm (i)} Conjecture \ref{pre-Bloch}  is true for $F=\Fbar$.
\par {\rm (ii)} Let $n \ge 3$. Then $\t_n(\Fbar)$ is a
torsion divisible group. In particular
\[
\ppp_n(\Fbar)\otimes \q=
\begin{cases}
\q & \text{if $n$ is odd}\\
0 & \text{if $n$ is even.}
\end{cases}
\]
\end{cor}
\begin{proof}
(i) By \ref{p234} we may assume $n \ge 4$.
The conjecture follows from Theorems \ref{friedlander-thm} and \ref{inj}.\\
(ii) Again we may assume $n \ge 4$.
Look at the spectral sequence $E_{r,s}^1(n)$ with $k=\q$. Since
${\overline{\F}}_q^\ast$ is torsion,
\begin{gather*}
\begin{array}{l}
E_{r, s}^1(n)= \begin{cases} H_s(\GL_{n-r}(\Fbar), \q) &
\text{if $0 \le r \le n$}\\
H_s(\GL_n(\Fbar), H_{n-1}(X_n, \q)) & \text{if $r=n+1$} \\
0 & \text{if $r \ge n+2.$}\end{cases}
\end{array}
\end{gather*}
Now by an easy analysis of this spectral sequence, using
\ref{friedlander-thm}, we see that $\ppp_n(\Fbar)$ is torsion.
The rest follows from
Lemma \ref{tn-pn}.
\end{proof}

\begin{cor}
Let $\char(F) \neq 0$.
Then the
following conditions are equivalent
\par {\rm (i)}
$\ppp_n(F)\otimes \z/l=
\begin{cases}
\z/l & \text{if $n$ is odd}\\
0 & \text{if $n$ is even}
\end{cases}$
for all $n$,
\par {\rm (ii)} $H_n(\GL_{n-1}, \z/l) \overset{\sim}{\larr}
H_n(\GL_{n}, \z/l)$ for all $n$.
\end{cor}
\begin{proof}
It is sufficient to prove that in part (i) of Theorem \ref{inj} the
surjectivity follows from the injectivity. If $\char(F)=p \neq 0$, then it
contains a copy of $\overline{\F}_p$. Consider the commutative diagram
\begin{gather*}
\begin{array}{ccc}
H_{n+1}(\GL_{n-1}(\overline{\F}_p), \z/l) & \larr &
H_{n+1}(\GL_{n-1}(F), \z/l) \\
\Big\downarrow & & \Big\downarrow \\
H_{n+1}(\GL_{n}(\overline{\F}_p), \z/l)&{\larr} &
H_{n+1}(\GL_{n}(F), \z/l)\\
\Big\downarrow & & \Big\downarrow \\
 H_{n+1}(\GL(\overline{\F}_p), \z/l)&\overset{\sim} {\larr} &
 H_{n+1}(\GL(F), \z/l).
\end{array}
\end{gather*}
By \ref{friedlander-thm}  the left column maps are bijective. By a
theorem of Suslin \cite[Cor. 1]{suslin1983} the bottom row map is
bijective. Now the claim follows easily.
\end{proof}

For a topological group $G$ let $BG^\top$ be its classifying space
with its underlying topology and $BG$ be its classifying space as
a topological group with discrete topology. By the functorial
property of $B$ we have a natural map $\psi: BG \arr BG^\top$.

\begin{conj}[Friedlander-Milnor Conjecture]\label{FMC}
Let $G$ be a Lie group. The canonical map
$\psi: BG \arr BG^\top$ induces isomorphism of homology and cohomology
with any finite abelian coefficient group.
\end{conj}

See \cite{milnor1983} and \cite{sah1989} for more information in
this direction. Here is a strong support for Conjecture \ref{pre-Bloch}.

\begin{prp}\label{FMC-pn}
The following conditions are equivalent
\par {\rm (i)}
$\ppp_n(\C)\otimes \z/l=
\begin{cases}
\z/l & \text{if $n$ is odd}\\
0 & \text{if $n$ is even}
\end{cases}$
for all $n$,
\par {\rm (ii)} $H_n(\GL_{n-1}(\C), \z/l) \overset{\sim}{\larr}
H_n(\GL_{n}(\C), \z/l)$ for all $n$,
\par {\rm (iii)}  $H_n(B\GL_{n-1}(\C), \z/l) \overset{\sim}{\larr}
H_n( B\GL_{n-1}(\C)^\top, \z/l)$ for all $n$.
\end{prp}
\begin{proof}
It is well-known that
\[
\GL_{n-1}(\C)\arr \GL_{n}(\C)\arr \GL_{n}(\C)/\GL_{n-1}(\C)
\]
is a fibration and $\GL_{n}(\C)/\GL_{n-1}(\C)$ is $(2n-2)$-connected
\cite[Thm. 3.15]{mimura-tota1991}. Hence
$\pi_i(\GL_{n-1}(\C)) \arr \pi_i(\GL_{n}(\C))$ is injective if
$i \le 2n-3$ and surjective if $i \le 2n-2$ which imply that
\[
\pi_j(B\GL_{n-1}(\C)^\top) \arr \pi_j(B\GL_{n}(\C)^\top)
\]
is injective if
$j \le 2n-2$ and surjective if $j \le 2n-1$. Therefore
\[
H_j(B\GL_{n-1}(\C)^\top, \z) \arr H_j(B\GL_{n}(\C)^\top, \z)
\]
is injective for $j \le 2n-2$ and surjective for $j \le 2n-1$
\cite[Chap. 7, Sec. 5. Thm. 9]{spanier1966}.
We call this topological stability.

Here we prove (i) $\Leftrightarrow$ (ii).
The proof of (ii) $\Leftrightarrow$ (iii) is similar and
easier.\\
(i) $\Rightarrow$ (ii) This immediately follows from \ref{inj}.\\
(ii ) $\Rightarrow$ (i) By \ref{inj} it is sufficient to prove that
\[
H_{n+1}(\GL_{n-1}(\C), \z/l) {\arr}
H_{n+1}(\GL_{n}(\C), \z/l)
\]
is surjective. For this we look at the following commutative diagram
\begin{gather*}
\begin{array}{ccc}
H_{n+1}(B\GL_{n-1}(\C), \z/l) & \larr &
H_{n+1}(B\GL_{n-1}(\C)^\top, \z/l) \\
\Big\downarrow & & \Big\downarrow \\
H_{n+1}(B\GL_{n}(\C), \z/l)&{\larr} &
H_{n+1}(B\GL_{n}(\C)^\top, \z/l)\\
\Big\downarrow & & \Big\downarrow \\
 H_{n+1}(B\GL(\C), \z/l)&\overset{\sim} {\larr} &
 H_{n+1}(B\GL(\C)^\top, \z/l).
\end{array}
\end{gather*}
By a theorem of Suslin the last row
is isomorphism \cite[Cor. 4.8]{suslin1984} and by topological
stability the column maps in the right are isomorphism for $n \ge
3$. By a result of Milnor, for any Lie group $G$ with a finite
number of connected components the map
\[
H_{i}(BG, \z/l) \arr  H_{i}(BG^\top, \z/l)
\]
is always surjective \cite[Thm. 1]{milnor1983}.
So the row maps of the diagram are surjective. By (ii) and \ref{stability}
the bottom column map in the left
is isomorphism,
so the middle row map is isomorphism. All these imply that
the first column map in the left of the diagram is surjective.
\end{proof}


\begin{rem}
The original goal of Loday to introduce the higher pre-Bloch groups
in \cite{loday1987} was that it might help one to study
$H_{n+1}(\GL_{n})/ H_{n+1}(\GL_{n-1})$,
which is motivated by the Bloch-Wigner exact sequence
and also by a result of Suslin which describes the quotient group
$H_{n}(\GL_{n})/ H_{n}(\GL_{n-1})$ explicitly \cite[Thm.
3.4]{suslin1985}.

It is easy to define a natural map
\[
\varrho_n: H_{n+1}(\GL_{n})/  H_{n+1}(\GL_{n-1}) \arr \ppp_n(F).
\]
This map can be constructed using exact
sequence (\ref{exact2}). From the short exact sequence
$0 \arr \partial_1(C_1(\fn)) \arr C_0(\fn) \arr \z \arr 0$
we get the connecting homomorphism
$H_{n+1}(\GL_n) \arr H_{n}(\GL_n,\partial_1(C_1(\fn)))$.
Iterating this process
we get a homomorphism $\eta_n: H_{n+1}(\GL_n) \arr \t_n(F)$.
Since the epimorphism $C_0(\fn) \arr \z \arr 0$ has a
$\GL_{n-1}$-equivariant section $m \mt m(\lan e_n \ran)$, the
restriction of $\eta_n$ to $H_{n+1}(\GL_{n-1})$ is zero. Thus we obtain
a homomorphism
\[
H_{n+1}(\GL_n)/H_{n+1}(\GL_{n-1}) \arr \t_n(F).
\]
The Composition of $\eta_n$ with the map $\t_n(F) \arr \ppp_n(F)$,
constructed in the previous section, gives us the map that we are
looking for. This map also can be constructed on the level of
complexes. For details of this approach see \cite{yagunov2000}.

In the light of Conjecture \ref{pre-Bloch},
it is convenient to ask the following question.

\begin{athm}[Question.]
Let $F$ be algebraically closed. Is $H_{n+1}(\GL_{n})/
H_{n+1}(\GL_{n-1})$ divisible?
\end{athm}
By \ref{friedlander-thm}, the answer to this question
is positive if $F=\Fbar$ and
in the next section we
show that the answer also is positive for $n \le 4$.
Using Theorem \ref{inj} one can show that Conjecture
\ref{pre-Bloch} gives a positive answer to the above question for
$n$ even
(see the proof of \ref{f-m-3}(ii)).
\end{rem}

\section{Lower degree homology groups}\label{lower-homology}

Here we demonstrate that the equivalence conditions in Theorem
\ref{inj} are true for $n \le 4$.
In this section we assume that $F$ is
algebraically closed, unless we mention it.

\begin{prp}\label{gli}
We have
\par
{\rm (i)}
$H_{2}(\GL_{1},\z/l) \overset{\sim}{\arr} H_{2}(\GL_{2},\z/l)$ and
$0=H_{3}(\GL_{1},\z/l) \two H_{3}(\GL_{2},\z/l)$
\par {\rm (ii)}
$H_{3}(\GL_{2},\z/l) \overset{\sim}{\arr} H_{3}(\GL_{3},\z/l)$ and
$H_{4}(\GL_{2},\z/l) \two H_{4}(\GL_{3},\z/l)$,
\par {\rm (iii)} $H_{4}(\GL_{3},\z/l)
\overset{\sim}{\arr} H_{4}(\GL_{4},\z/l)$ and
$H_{5}(\GL_{3},\z/l) \two H_{5}(\GL_{4},\z/l)$,
\par {\rm (iv)} $H_{4}(\GL_{3})/ H_{4}(\GL_{2})$ and
$H_{5}(\GL_{4})/ H_{5}(\GL_{3})$ are divisible.
\end{prp}
\begin{proof}
By Theorem \ref{inj}, to proof (i), (ii) and (iii) it is sufficient to prove
that $\ppp_2(F)$ and $\ppp_4(F)$ are divisible and
$\ppp_3(F)\otimes\z/l=\z/l$.
Dupont and Sah \cite[Thm. 5.1]{dupont-sah1982} proved that
$\ppp_2(F)$ is divisible. The rest follows from this and Lemma \ref{p234}.
We should mention that to prove Lemma \ref{p234} we used the main
results of \cite{mirzaii2004} and \cite{mirzaii2005}, which are
difficult since are on an arbitrary infinite field. The proof of those
results has great simplification over an algebraically closed
field and homology with $\z/l$ coefficient, for example the proof
of Lemmas 5.2, 5.3 and 5.4 in \cite{mirzaii2005} are easy (some
even trivial) as $H_{2i+1}(\fff, \z/l)=0$ for $i \ge 0$.
\\
(iv) Consider the following commutative diagram with exact rows
\begin{gather*}
\begin{array}{ccccccccc}
0\!\! &\!\! \arr H_4(\GL_2)\otimes \z/l & \!\!\arr & \!\!
 H_4(\GL_2, \z/l) &\!\! \arr
& \!\! \tors(H_3(\GL_2), \z/l)& \arr & \!\!0  \\
& \Big\downarrow & & \Big\downarrow & & \Big\downarrow & \\
0\!\! &\!\! \arr H_4(\GL_3)\otimes \z/l &\!\! \arr & \!\!
 H_4(\GL_3, \z/l) &\!\! \arr
& \!\! \tors(H_3(\GL_3), \z/l)&\!\! \arr & \!\! 0.
\end{array}
\end{gather*}
Since $H_3(\GL_3)\simeq H_3(\GL_2)\oplus K_3^M(F)$ \cite[Cor.
5.5]{mirzaii2004}, and since  $K_i^M(F)$ is uniquely divisible for
$i\ge 2$ \cite[1.2]{bass-tate1973}, the right column map is
isomorphism. The middle column map is surjective by (ii). Thus the
left column map is surjective too. This shows that
$H_{4}(\GL_{3})/ H_{4}(\GL_{2})$ is $l$-divisible. The proof of
divisibility of the group $H_{5}(\GL_{4})/H_{5}(\GL_{3})$ is
analogue and for this one should use (iii) and \cite[Cor.5.7]{mirzaii2005}.
\end{proof}
~
\begin{prp}\label{sli}
{\rm (i)} $H_{3}(\SL_{2},\z/l)=H_{3}(\SL_{3},\z/l)=0$,
so $H_{3}(\SL_{2})$ and $H_{3}(\SL_{3})$ are divisible,
\par {\rm (ii)} $H_{4}(\SL_{2},\z/l) \two H_{4}(\SL_{3},\z/l)
\overset{\sim}{\arr} H_{4}(\SL_{4},\z/l)$ and
$H_{5}(\SL_{3},\z/l) \two H_{5}(\SL_{4},\z/l)$,
\par {\rm (iii)} $H_{4}(\SL_{3})/ H_{4}(\SL_{2})$ and
$H_{5}(\SL_{4})/ H_{5}(\SL_{3})$ are divisible.
\end{prp}
\begin{proof}
Since $H_i(\SL, \z/l) \harr H_i(\GL, \z/l)$ for all i,
by homology stability theorem \ref{stability}
and Prop. \ref{gli}(i) one gets $H_3(\SL_3, \z/l)=0$.
The exact sequence
\[
1 \arr \SL_n \arr \GL_n \arr \fff \arr 1
\]
induces the Lyndon-Hochschild-Serre spectral sequence
\[
{}_nE_{p, q}^2=H_p(\fff, H_q(\SL_n, \z/l)) \Rightarrow
H_{p+q}(\GL_n, \z/l).
\]
It is easy to see that $H_q(\SL_n, \z/l)=0$ for $q=1,2$, thus
${}_nE_{p, q}^2=0$ for $q=1,2$. Triviality of
$H_i(\fff,\z/l)$ for $i$ odd, implies that for $p$ odd
\[
{}_nE_{p, q}^2=H_p(\fff, H_q(\SL_n, \z/l))=
H_p(\fff,\z/l)\otimes H_q(\SL_n, \z/l)=0.
\]
Since the above exact sequence splits,
${}_nd_{p, 0}^r$ are trivial maps.
Thus ${}_2E_{0, 3}^\infty={}_2E_{0, 3}^2= H_3(\SL_2,
\z/l)=H_3(\GL_2, \z/l)=0$, which imply ${}_2E_{2, 3}^2=0$. So we
obtain the exact sequences
\begin{equation}\label{exact7}
0 \arr H_{4}(\SL_{2},\z/l) \arr
H_{4}(\GL_{2},\z/l) \arr H_{4}(\fff,\z/l) \arr 0,
\end{equation}
\begin{equation}\label{exact8}
{}_2E_{2, 4}^2 \arr H_{5}(\SL_{2},\z/l) \arr
H_{5}(\GL_{2},\z/l) \arr 0.
\end{equation}
With an analogue argument for ${}_nE_{p, q}^2$, $n=3, 4$,
we obtain the exact sequences
\begin{equation}\label{exact9}
0 \arr H_{4}(\SL_{3},\z/l) \arr
H_{4}(\GL_{3},\z/l) \arr H_{4}(\fff,\z/l) \arr 0,
\end{equation}
\begin{equation}\label{exact10}
{}_3E_{2, 4}^2 \arr H_{5}(\SL_{3},\z/l) \arr
H_{5}(\GL_{3},\z/l) \arr 0,
\end{equation}
\begin{equation}\label{exact11}
0 \arr H_{4}(\SL_{4},\z/l) \arr
H_{4}(\GL_{4},\z/l) \arr H_{4}(\fff,\z/l) \arr 0,
\end{equation}
\begin{equation}\label{exact12}
{}_4E_{2, 4}^2 \arr H_{5}(\SL_{4},\z/l) \arr
H_{5}(\GL_{4},\z/l) \arr 0.
\end{equation}
If $m \le n$, then there is a natural map of spectral sequences
\[
{}_mE_{p, q}^2 \arr {}_nE_{p, q}^2.
\]
Now the isomorphism $H_{4}(\SL_{3},\z/l)
\overset{\sim}{\arr} H_{4}(\SL_{4},\z/l)$ can be deduced from the natural
map from exact sequence (\ref{exact9}) to exact sequence  (\ref{exact11})
and the corresponding result for $\GL$ in Proposition \ref{gli}.
This isomorphism
implies that ${}_3E_{2, 4}^2 \overset{\sim}{\arr} {}_4E_{2, 4}^2$.
From this, \ref{gli}(iii) and exact sequences (\ref{exact10})
and  (\ref{exact12}) we obtain the surjectivity
$H_{5}(\SL_{3},\z/l) \two H_{5}(\SL_{4},\z/l)$.
The proof of (iii) is analogue to the case $\GL$ in \ref{gli}(iv)
using \cite[Cor. 6.2]{mirzaii2004} and \cite[Prop. 5.8]{mirzaii2005}.
\end{proof}

Here is a new case of the Friedlander-Milnor Conjecture.

\begin{cor}\label{f-m-2}
Let $G=\GL_3(\C)$ or $\SL_3(\C)$. Then for any finite abelian groups $A$,
\[
H_4(BG, A)\overset{\sim}{\larr} H_4(BG^\top, A).
\]
\end{cor}
\begin{proof}
We may assume $A=\z/l$, where $l$ is a prime. Now a
similar argument as in the proof of Proposition
\ref{FMC-pn} using \ref{gli} and  \ref{sli} will prove this claim.
\end{proof}

\begin{cor}\label{f-m-3}
Let $n \ge 3$. Then
\par {\rm (i)}
$H_4(\SL_n, \z/l)\simeq
\begin{cases}
0 & \text{if $\char(F)| l$}\\
\z/l &\text{if $\char(F) \nmid$ l,}
\end{cases}$
\par {\rm (ii)} $H_4(\SL_n)$ is uniquely divisible.
\end{cor}
\begin{proof}
We may assume that $l$ is a prime. By a result of Suslin, the
$K$-theory of algebraically closed fields with finite
coefficient, $K_i(F,\z/l)$, does not depend on the field and
$K_i(F,\z/l)$ is trivial if
\par (1)  $i$ is odd,
\par (2)  $i \ge 1$ when $\char(F)=l \neq 0$\\
(see \cite{suslin1983} and  \cite[Cor. 3.13]{suslin1984}). This
implies that the group $H_i(\SL, \z/l)$ does not depends on $F$
and $H_i(\SL, \z/l)$ is trivial in the above cases
(see \cite[Cor. 1, Cor. 2]{suslin1983}).\\
(i) To prove this claim
it is sufficient to prove it for $F=\C$. It is well-known that
$\SL_n(\C)$, as a Lie group, is $2$-connected and
$\pi_3(\SL_n(\C)) \simeq \z$. This implies that $B\SL_n(\C)^\top$
is 3-connected and $\pi_4(B\SL_n(\C)^\top) \simeq \z$. Therefore
\[
H_4(B\SL_n(\C)^\top, \z) \simeq \z.
\]
From  Cor. \ref{f-m-2} we have
\[
H_4(\SL_n(\C), \z/l)\simeq H_4(B\SL_n(\C), \z/l)\simeq
H_4(B\SL_n(\C)^\top, \z/l)\simeq \z/l.
\]
(ii) The exact sequence
$0\arr \z \overset {l.}{\arr} \z \arr \z/l \arr 0$ induces the long exact
sequence
\begin{gather*}
\cdots \arr H_4(\SL)\overset {l.}{\arr}
H_4(\SL) \arr H_4(\SL, \z/l)
\arr H_3(\SL)
\overset {l.}{\arr}H_3(\SL) \arr \cdots.
\end{gather*}
For $n \ge 4$ the claim follows from (i), the triviality of
$H_{5}(\SL, \z/l)$ (and of $H_{4}(\SL, \z/l)$ if $\char(F)=l \neq 0$) and
the following fact
\[
H_3(\SL)\simeq V \oplus
\begin{cases}  \q/\z&
\text{if $\char(F)=0$}  \\
 \q/\z[\frac{1}{p}] &
\text{if $\char(F)=p \neq 0$,}
\end{cases}
\]
where $V$ is a uniquely divisible group \cite{suslin1991}. The case $n=3$
follows from this and the fact that
$H_4(\SL_4) \simeq H_4(\SL_3) \oplus K_4^M(F)$
\cite[Prop. 5.8]{mirzaii2005}. Note that $K_4^M(F)$ is uniquely divisible.
\end{proof}

\begin{exa}\label{pn-R}
(i) Propositions \ref{gli} and \ref{sli} are true if $F=\R$ and $2
\nmid l$. Because then $\ppp_2(\R)$ is divisible \cite[2.14, 4.1(a)]{sah1989}
and $H_{i}(\R^\ast, \z/l)=0$ for $i \ge 1$. For example, a similar
argument as in the above shows that the groups
\[
H_4(\GL_3(\R))/H_4(\GL_2(\R)) \ \
{\rm and} \ \ H_5(\GL_4(\R))/ H_5(\GL_3(\R))
\]
are $l$-divisible (see \cite[Cor. 5.5]{mirzaii2004}, \cite[Example 1]{mirzaii2005}).
The same is true if one replaces $\GL$ with $\SL$.
%
\par (ii) Corollary \ref{f-m-3}, homology stability theorem \ref{stability}
and exact sequences (\ref{exact9}) and (\ref{exact11}) imply that
for $n \ge 3$
\[
H_4(\GL_n, \z/l)\simeq \z/l \oplus H_4(\fff, \z/l)
\simeq \z/l \oplus \z/l.
\]
Using \ref{f-m-3} and \cite[Cor. 5.7]{mirzaii2005} it is easy to
prove that $H_4(\GL_n)$ is uniquely divisible for $n \ge 3$.
\end{exa}

\section{Some generalizations}\label{generalization}

One can generalize Theorem \ref{inj} as follows;

\begin{prp}\label{inj2}
Let F be an infinite field and let $k=\q$ or $k=\z/l$, $l$ a
prime, such that $K_2^M(F) \otimes k=0$.
Then the
following are equivalent;
\par {\rm (i)} For $3 \le m \le n$ the map
$H_{m}(\GL_{m-1}, k) {\arr} H_{m}(\GL_{m},k)$
is injective and
the map
$H_{m+1}(\GL_{m-1}, k) {\arr} H_{m+1}(\GL_{m},k)$ is
surjective,
\par {\rm (ii)} For $3 \le m \le n$, the complex
$\t_{m-2}^{(k)}(F) \arr \t_{m}^{(k)}(F) \arr H_1(\fff, k)
\arr 0$ is exact.
\end{prp}
\begin{proof}
Clearly  $K_i^M(F) \otimes k=0$ for $i \ge 2$.
The proof is similar to the proof of Theorem \ref{inj}. We leave the
details to the reader.
\end{proof}

\begin{exa}\label{F-k}
Here are some examples of pairs of fields $(F, k)$ such that
$K_2^M(F) \otimes  k=0$:
\par (1) $F$ any global field and $k=\q$ (see \cite{bass-tate1973}),
\par (2) $F$ a perfect field of $\char(F)=l \neq 0$ and $k=\z/l$
(see \cite{bass-tate1973}),
\par (3) $F$ algebraically closed and $k=\z/l$
(see Sections \ref{finite-coeff} and \ref{lower-homology}),
\par (4) $F$ a local field and  $k=\z/l$, $ l \nmid  |\mu(F)|$
(see \cite[Example 1.7]{milnor1970}),
\par (5) $F=\R$ and $k=\z/l$, $l \neq 2$
(see \cite[Example 1.6]{milnor1970}),
\par (6) $F=\overline{\F}_q$ and $k$ any prime field.
(see \cite[Example 1.5]{milnor1970})
\par (7) $F=\overline{\q}$ and $k$ any prime field
\cite[Cor. 10.21]{dupont2001}.

~

To give more examples first we state a result of Milnor \cite[\S
2]{milnor1970}.

\begin{athm}[Theorem.]
For every field $F$ we have an exact sequence
\[
0 \arr K_{n+1}^M(F) \arr K_{n+1}^M(F(t)) \arr
\coprod_{P \in \Spec(F[t])} K_{n}^M(F[t]/P) \arr 0.
\]
\end{athm}

The following pairs of fields $(F, k)$ with the desired property
follow from this theorem and the above cases;

~

\par (8) $F=\overline{\F}_q(T)$ or $F={\F}_q(T)$
and $k=\q$ or $k=\z/l$, $l \neq \char(F)$,
\par (9) $F=E(t)$, $E$ algebraically closed and $k=\z/l$,
\par (10) $F=\R(t)$ and $k=\z/l$, $l\neq 2$.
\end{exa}

By Lemma \ref{tn-pn}, Conjecture \ref{pre-Bloch} immediately
follows from the following conjecture.

\begin{conj}\label{H(X_N)}
Let $n \ge 3$. Then
$\t_{n-2}^{(k)}(F) \arr \t_{n}^{(k)}(F) \arr H_1(\fff, k) \arr 0$
is exact.
\end{conj}

By Lemma \ref{p234} this conjecture is true for $n=3,4$. By
\ref{pn-fbar} it is also true for the algebraic closure of a
finite field when $k$ is a prime field. We should mention that the
surjectivity of $\t_{n}^{(k)}(F) \arr H_1(\fff, k)$ is proven in
Section \ref{Pre-Bloch}.
We have the following results analogue to Prop. \ref{gli}.



\begin{prp}\label{F-k-gl}
Let the pair $(F, k)$ be as in Example \ref{F-k}. Then
\par {\rm (i)}
$H_{3}(\GL_{2},k) \overset{\sim}{\arr} H_{3}(\GL_{3},k)$ and
$H_{4}(\GL_{2},k) \two H_{4}(\GL_{3},k)$,
\par {\rm (ii)} $H_{4}(\GL_{3},k)
\overset{\sim}{\arr} H_{4}(\GL_{4},k)$ and
$H_{5}(\GL_{3},k) \two H_{5}(\GL_{4},k)$.
\end{prp}
\begin{proof}
The proof is similar to the proof of \ref{gli} using \ref{p234} and \ref{inj2}.
\end{proof}
~
\\
~

\subsection*{Acknowledgement}
This article was written during my stay at the Mathematics
Department of Bielefeld University and Queen's University Belfast.
I would like to thank them for their support and hospitality. I
would especially like to thank professor A. Bak for his
encouragement and his interest in this work. Finally I would like to
thank W. van der Kallen for his interest in this work.


\bigskip

\address{{\footnotesize
\
\\
Department of Pure Mathematics\\
Queen's University\\
Belfast  BT7 1NN \\
Northern Ireland\\
Email:\ bmirzaii@gmail.com
}}

\begin{thebibliography}{99}



\bibitem{bass-tate1973} Bass, H., Tate, J.
The Milnor ring of a global field. 349--446. Lecture Notes in Math.,
Vol. {\bf 342}, 1973, 349--446.




\bibitem{dupont2001} Dupont, J- L. Scissors congruences, group homology and
characteristic classes. Nankai Tracts in Mathematics, 1,
World Scientific Publishing Co., Inc., River Edge, NJ, 2001.

\bibitem{dupont-sah1982} Dupont, J- L.; Sah, C. Scissors congruences. II.
J. Pure Appl. Algebra  {\bf 25}  (1982), no. 2, 159--195.




\bibitem{friedlander1976} Friedlander, E. M. Homological stability for
classical groups over finite fields.  Algebraic $K$-theory,
 pp. 290--302. Lecture Notes in Math., Vol. {\bf 551},
 1976.




\bibitem{loday1987}
Loday, J. L. Comparaison des homologies du groupe lin\'eaire
et de son alg\`ebre de Lie.
Ann. Inst. Fourier (Grenoble)  {\bf 37}  (1987),  no. 4, 167--190.


\bibitem{milnor1970}
Milnor, J. Algebraic $K$-theory and quadratic forms.  Invent. Math.
{\bf 9} 1970, 318--344.

\bibitem{milnor1983} Milnor, J. On the homology of Lie groups made discrete.
Comment. Math. Helv. {\bf 58} (1983), no. 1, 72--85.

\bibitem{mimura-tota1991} Mimura, M., Toda, H. Topology of Lie groups. I, II.
Translations of Mathematical
Monographs, 91. American Mathematical Society, Providence, RI, 1991.

\bibitem{mirzaii2003} Mirzaii, B. Homology stability for unitary
groups II. $K$-Theory {\bf 36} (2005), no. 3--4, 305--326.

\bibitem{mirzaii2004} Mirzaii, B. Third homology of general linear groups,
Preprint. To appear.

\bibitem{mirzaii2005} Mirzaii, B. Homology of $\SL_n$ and $\GL_n$ over an
infinite field. Preprint, available at http://arxiv.org/abs/math/0605722


\bibitem{sah1989} Sah, C. Homology of classical Lie groups made discrete.
III. J. Pure Appl. Algebra {\bf 56} (1989), no. 3, 269--312.

\bibitem{spanier1966} Spanier E. H, Algebraic Topology. McGram Hill (1966).

\bibitem{suslin1983} Suslin, A. A. On the $K$-theory of algebraically
closed fields.  Invent. Math.  {\bf 73}  (1983),  no. 2, 241--245.

\bibitem{suslin1984} Suslin, A. A. On the $K$-theory of local fields.
J. Pure Appl. Algebra {\bf 34}  (1984),  no. 2-3, 301--318.

\bibitem{suslin1985} Suslin, A. A. Homology of ${\rm GL}\sb{n}$,
characteristic classes and Milnor $K$-theory. Proc. Steklov Math.
{\bf 3} (1985), 207--225.

\bibitem{suslin1991} Suslin A. A. $K\sb 3$ of a field, and the Bloch group.
Proc. Steklov Inst. Math. 1991, {\bf 183} no. 4, 217--239.



\bibitem{yagunov2000} Yagunov, S. On the homology of $\GL_n$ and the higher
pre-Bloch groups. Canad J. Math. {\bf 52} (6), 2000, 1310--1338.

\end{thebibliography}
\end{document}